# DOUBLE CONVERGENCE OF A FAMILY OF DISCRETE DISTRIBUTED MIXED ELLIPTIC OPTIMAL CONTROL PROBLEMS WITH A PARAMETER

*Domingo A. Tarzia*

CONICET - Depto. Matemática,
FCE, Univ. Austral, Paraguay 1950,
S2000FZF Rosario, Argentina.
Tel.: +54-341-5223093; Fax: +54-341-5223001
E-mail: DTarzia@austral.edu.ar

**ABSTRACT**

We consider a bounded domain $\Omega \subset \mathbb{R}^n$ whose regular boundary $\Gamma = \partial\Omega = \Gamma_1 \cup \Gamma_2$ consists of the union of two disjoint portions $\Gamma_1$ and $\Gamma_2$ with $\operatorname{meas}(\Gamma_1) > 0$. The convergence of a family of continuous distributed mixed elliptic optimal control problems ($P_\alpha$), governed by elliptic variational equalities, when the parameter $\alpha$ of the family (the heat transfer coefficient on the portion of the boundary $\Gamma_1$) goes to infinity was studied in Gariboldi - Tarzia, Appl. Math. Optim., 47 (2003), 213-230. It has been proved that the optimal control, and their corresponding system and adjoint system states are strongly convergent, in adequate functional spaces, to the optimal control, and the system and adjoint states of another distributed mixed elliptic optimal control problem ($P$) governed also by an elliptic variational equality with a different boundary condition on the portion of the boundary $\Gamma_1$

We consider the discrete approximations ($P_{h\alpha}$) and ($P_h$) of the optimal control problems ($P_\alpha$) and ($P$) respectively, for each $h > 0$ and for each $\alpha > 0$, through the finite element method with Lagrange's triangles of type 1 with parameter $h$ (the longest side of the triangles). We also discretize the elliptic variational equalities which define the system and their adjoint system states, and the corresponding cost functional of the distributed optimal control problems ($P_\alpha$) and ($P$). The goal of this paper is to study the double convergence of this family of discrete distributed mixed elliptic optimal control problems ($P_{h\alpha}$) when the parameters $\alpha$ goes to infinity and the parameter $h$ goes to zero simultaneously. We prove the convergence of the discrete optimal controls, the discrete system and adjoint states of the family ($P_{h\alpha}$) to the corresponding to the discrete distributed mixed elliptic optimal control problem ($P_h$) when $\alpha \to \infty$, for each $h > 0$, in adequate functional spaces. We study the convergence of the discrete distributed optimal control problems ($P_{h\alpha}$) and ($P_h$) when $h \to 0$ obtaining a commutative diagram which relates the continuous and discrete distributed mixed elliptic optimal control problems $(P_{h\alpha}), (P_\alpha), (P_h)$ and ($P$) by taking the limits $h \to 0$ and $\alpha \to \infty$ respectively. We also study the double convergence of ($P_{h\alpha}$) to ($P$) when $(h,\alpha) \to (0, +\infty)$ which represents the diagonal convergence in the above commutative diagram.

I. **Introduction**

The purpose of this paper is to do the numerical analysis, by using the finite element method, of the convergence of the continuous distributed mixed optimal control problems with respect to a parameter (the heat transfer coefficient) given in [10, 11] obtaining a double convergence when the parameter of the finite element method goes to zero and the heat transfer coefficient goes to infinity.

We consider a bounded domain $\Omega \subset \mathbb{R}^n$ whose regular boundary $\Gamma = \partial\Omega = \Gamma_1 \cup \Gamma_2$ consists of the union of two disjoint portions $\Gamma_1$ and $\Gamma_2$ with $\operatorname{meas}(\Gamma_1) > 0$. We consider the following elliptic partial differential problems with mixed boundary conditions, given by:

$$-\Delta u = g \quad \text{in } \Omega \; ; \quad u = b \quad \text{on } \Gamma_1 \; ; \quad -\frac{\partial u}{\partial n} = q \quad \text{on } \Gamma_2 \tag{1}$$

and

$$-\Delta u = g \quad \text{in } \Omega ; \quad -\frac{\partial u}{\partial n} = \alpha(u-b) \quad \text{on } \Gamma_1 ; \quad -\frac{\partial u}{\partial n} = q \quad \text{on } \Gamma_2 \qquad (2)$$

where $g$ is the internal energy in $\Omega$, $b = Const. > 0$ is the temperature on $\Gamma_1$ for the system (1) and the temperature of the external neighborhood on $\Gamma_1$ for the system (2) respectively, $q$ is the heat flux on $\Gamma_2$ and $\alpha > 0$ is the heat transfer coefficient on $\Gamma_1$. The systems (1) and (2) can represent the steady-state two-phase Stefan problem for adequate data [21, 22]. We consider the following continuous distributed optimal control problem ($P$) and a family of continuous distributed optimal control problems ($P_\alpha$) for each parameter $\alpha > 0$, defined in [10], where the control variable is the internal energy $g$ in $\Omega$, that is: Find the continuous distributed optimal controls $g_{op} \in H = L^2(\Omega)$ and $g_{\alpha_{op}} \in H$ (for each $\alpha > 0$) such that:

$$\text{Problem (P): } J(g_{op}) = \min_{g \in H} J(g), \quad \text{Problem } (P_\alpha): J_\alpha(g_{\alpha_{op}}) = \min_{g \in H} J_\alpha(g), \qquad (3)$$

where the quadratic cost functional $J : H \to \mathbb{R}_0^+$ and $J_\alpha : H \to \mathbb{R}_0^+$ are defined by the following expresions [2, 18, 26]:

$$a) J(g) = \frac{1}{2}\|u_g - z_d\|_H^2 + \frac{M}{2}\|g\|_H^2, \quad b) J_\alpha(g) = \frac{1}{2}\|u_{\alpha g} - z_d\|_H^2 + \frac{M}{2}\|g\|_H^2 \qquad (4)$$

with $M > 0$ and $z_d \in H$ given, $u_g \in K$ and $u_{\alpha g} \in V$ are the state of the systems defined by the mixed ellliptic differential problems (1) and (2) respectively whose elliptic variational equality are given by [16]:

$$u_g \in K: \quad a(u_g, v) = (g, v)_H - \int_{\Gamma_2} qv d\gamma, \quad \forall v \in V_0 \qquad (5)$$

$$u_{\alpha g} \in V: \quad a_\alpha(u_{\alpha g}, v) = (g, v)_H - \int_{\Gamma_2} qv d\gamma + \alpha \int_{\Gamma_1} bv d\gamma, \quad \forall v \in V \qquad (6)$$

and their adjoint system states $p_g \in V$ and $p_{\alpha g} \in V$ are defined by the following elliptic variational equalities:

$$p_g \in V_o: \quad a(p_g, v) = (u_g - z_d, v), \quad \forall v \in V_0 \qquad (7)$$

$$p_{\alpha g} \in V: \quad a_\alpha(p_{\alpha g}, v) = (u_{\alpha g} - z_d, v), \quad \forall v \in V \qquad (8)$$

with:

$$V = H^1(\Omega), V_0 = \{v \in V, v/\Gamma_1 = 0\}, K = \{v \in V, v/\Gamma_1 = b\} = b + V_0, H = L^2(\Omega), Q = L^2(\Gamma_2) \qquad (9)$$

$$a(u, v) = \int_\Omega \nabla u . \nabla v dx, \quad a_\alpha(u, v) = a(u, v) + \alpha \int_{\Gamma_1} uv d\gamma, \quad (u, v) = \int_\Omega uv \, dx \qquad (10)$$

where the bilinear, continuous and symmetric forms $a$ and $a_\alpha$ are coercive on $V_0$ and $V$ respectively, that is [16]:

$$\exists \lambda > 0 \text{ such that } \lambda \|v\|_V^2 \leq a(v, v), \quad \forall v \in V_0. \qquad (11)$$

$$\exists \lambda_\alpha = \lambda_1 \min(1, \alpha) > 0 \text{ such that } \lambda_\alpha \|v\|_V^2 \leq a_\alpha(v, v), \quad \forall v \in V, \qquad (12)$$

and $\lambda_1 > 0$ is the coercive constant for the bilinear form $a_1$ [16, 21].

The unique continuous distributed optimal energies $g_{op}$ and $g_{\alpha_{op}}$ have been characterized in [10] as a fixed point on $H$ for a suitable operators $W$ and $W_\alpha$ over their optimal adjoint system states $p_{g_{op}} \in V_0$ and $p_{\alpha g_{\alpha_{op}}} \in V$ defined by:

$$W : H \to H \quad \text{such that} \quad W(g) = -\frac{1}{M} p_g, \qquad (13)$$

$$W_\alpha : H \to H \quad \text{such that} \quad W_\alpha(g) = -\frac{1}{M} p_{\alpha g}. \qquad (14)$$

The limit of the optimal control problem ($P_\alpha$) when $\alpha \to \infty$ was studied in [10] and it was proven that:

$$\lim_{\alpha \to \infty} \|u_{\alpha g_{\alpha_{op}}} - u_{g_{op}}\|_V = 0, \quad \lim_{\alpha \to \infty} \|p_{\alpha g_{\alpha_{op}}} - p_{g_{op}}\|_V = 0, \quad \lim_{\alpha \to \infty} \|g_{\alpha_{op}} - g_{op}\|_H = 0 \qquad (15)$$

for a large constant $M > 0$ by using the characterization of the optimal controls through fixed points (13) and (14); this restrictive hypothesis on data was eliminated in [11] by using the variational formulations. We can summary the conditions (15) saying that the distributed optimal control problems ($P_\alpha$) converges to the distributed optimal control problem (P) when $\alpha \to +\infty$.

Now, we consider the finite element method and a polygonal domain $\Omega \subset \mathbb{R}^n$ with a regular triangulation with Lagrange triangles of type 1, constituted by affine-equivalent finite element of class $C^0$ being $h$ the parameter of the finite element approximation which goes to zero [3,7]. Then, we discretize the elliptic variational equalities for the system states (6) and (5), the adjoint system states (8) and (7), and the cost functional (4a,b) respectively. In general, the solution of a mixed elliptic boundary problem belongs to $H^r(\Omega)$ with $1 < r \leq 3/2 - \varepsilon$ ($\varepsilon > 0$) but there exist some examples which solutions belong to $H^r(\Omega)$ with $2 \leq r$ [1, 17, 20]. Note that mixed boundary conditions play an important role in various applications, e.g. heat conduction and electric potential problems [12].

The goal of this paper is to study the numerical analysis, by using the finite element method, of the convergence results (15) corresponding to the continuous distributed elliptic optimal control problems $(P_\alpha)$ and $(P)$ when $\alpha \to +\infty$. The main result of this paper can be characterized by the following result:

**Theorem 1**

*The following commutative diagram which relates the continuous distributed mixed optimal control problems $(P_\alpha)$ and $(P)$, with the discrete distributed mixed optimal control problems $(P_{h\alpha})$ and $(P_h)$ is obtained by taking the limits $h \to 0$, $\alpha \to +\infty$ and $(h,\alpha) \to (0,+\infty)$ as follows:*

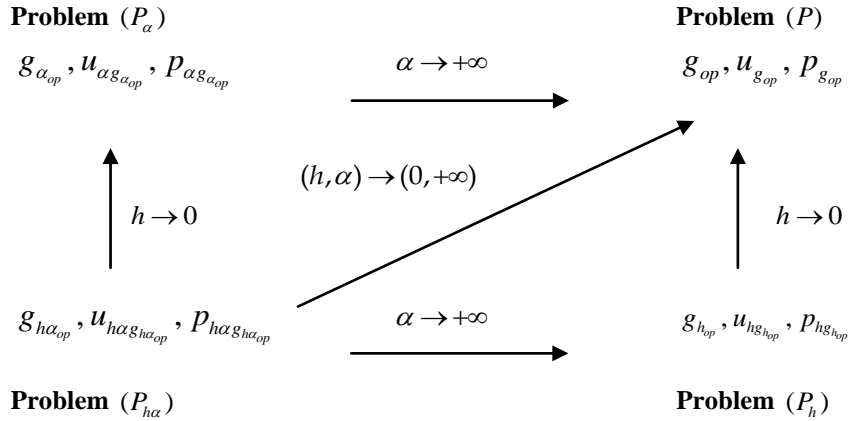

*where $u_{h\alpha g_{h\alpha_{op}}}$ and $p_{h\alpha g_{h\alpha_{op}}}$ are respectively the system and the adjoint system states of the discrete distributed mixed optimal control problem $(P_{h\alpha})$ for each $h > 0$ and $\alpha > 0$. Moreover, we obtain error estimates for the convergence when $h \to 0$ between the solution of problem $(P_{h\alpha})$ with respect to problem $(P_\alpha)$ for each $\alpha > 0$, and between the solution of problem $(P_h)$ with respect to problem $(P)$ respectively by using the fixed point characterization (13) and (14) of the optimal control problems $(P_h)$ and $(P_{h\alpha})$.*

The study of the limit $h \to 0$ of the discrete solutions of optimal control problems can be considered as a classical limit, see [4-6, 8, 9, 13-15, 19, 23, 24, 27, 28] but the limit $\alpha \to +\infty$, for each $h > 0$, and the double limit $(h,\alpha) \to (0,+\infty)$ can be considered as a new ones.

The paper is organized as follows. In Section II we define the discrete elliptic variational equalities for the state systems $u_{hg}$ and $u_{h\alpha g}$, we define the discrete distributed cost functional $J_h$ and $J_{h\alpha}$, we define the discrete distributed optimal control problems $(P_h)$ and $(P_{h\alpha})$, and the discrete elliptic variational equalities for the adjoint state systems $p_{hg}$ and $p_{h\alpha g}$ for each $h > 0$ and $\alpha > 0$. We obtain properties for the optimal control problem $(P_h)$: for system state $u_{hg}$ and adjoint system state $p_{hg}$, for the discrete cost functional $J_h$ and its corresponding optimality condition. We also define a contraction operator $W_h$ which allows obtain the optimal control $g_{h_{op}}$ as a fixed point.

We also obtain properties for the optimal control problem $(P_{h\alpha})$: for system $u_{h\alpha g}$ and adjoint system states $p_{h\alpha g}$, for the discrete cost functional $J_{h\alpha}$ and its corresponding optimality condition. We also define a contraction operator $W_{h\alpha}$ which allows obtain the optimal control $g_{h\alpha_{op}}$ as a fixed point.

In Section III we study the classical convergence of the discrete elliptic variational equalities for $u_{hg}, u_{h\alpha g}, p_{hg}$, and $p_{h\alpha g}$ as $h \to 0$ when $g$ is fixed (for each $\alpha > 0$). We study the convergences of the discrete distributed optimal control problems $(P_h)$ to $(P)$, and $(P_{h\alpha})$ to $(P_\alpha)$ when $h \to 0$ (for each $\alpha > 0$). We also study the error estimates for the optimal control problems $(P_h)$ and $(P_{h\alpha})$ (for each $\alpha > 0$) and the estimations for the discrete cost functional $J_h$ and $J_{h\alpha}$ corresponding to the discrete distributed optimal control problems $(P_h)$ and $(P_{h\alpha})$, for each $\alpha > 0$.

In Section IV we study the new convergence of the discrete distributed optimal control problems $(P_{h\alpha})$ to $(P_h)$ when $\alpha \to +\infty$ for each $h > 0$ and we obtain a commutative diagram which relates the continuous and discrete distributed mixed optimal control problems $(P_{h\alpha}), (P_\alpha), (P_h)$ and $(P)$ by taking the limits $h \to 0$ and $\alpha \to +\infty$.

In Section V we study the new double convergence of the discrete distributed optimal control problems $(P_{h\alpha})$ to $(P)$ when $(h, \alpha) \to (0, +\infty)$ and we obtain the diagonal convergence in the previous commutative diagram.

## II. Discretization by Finite Element Method and Properties

We consider the finite element method and a polygonal domain $\Omega \subset \mathbb{R}^n$ with a regular triangulation with Lagrange triangles of type 1, constituted by affine-equivalent finite element of class $C^0$ being $h$ the parameter of the finite element approximation which goes to zero [3, 7]. We can take $h$ equal to the longest side of the triangles $T \in \tau_h$ and we can approximate the sets $V, V_0$ and $K$ by:

$$V_h = \left\{ v_h \in C^0(\overline{\Omega}) / v_h/T \in P_1(T), \forall T \in \tau_h \right\}, \quad V_{0h} = \{ v_h \in V_h / v_h/\Gamma_1 = 0 \}; \quad K_h = b + V_{0h} \qquad (16)$$

where $P_1$ is the set of the polynomials of degree less than or equal to 1. Let $\pi_h : C^0(\overline{\Omega}) \to V_h$ be the corresponding linear interpolation operator. Then there exists a constant $c_0 > 0$ (independent of the parameter $h$) such that [3]:

$$a) \ \|v - \pi_h(v)\|_H \leq c_0 h^r \|v\|_r, \quad b) \ \|v - \pi_h(v)\|_V \leq c_0 h^{r-1} \|v\|_r, \quad \forall v \in H^r(\Omega), \ 1 < r \leq 2. \qquad (17)$$

We define the discrete cost functional $J_h, J_{h\alpha} : H \to \mathbb{R}_0^+$ by the following expressions:

$$a) \ J_h(g) = \frac{1}{2} \|u_{hg} - z_d\|_H^2 + \frac{M}{2} \|g\|_H^2, \quad b) \ J_{h\alpha}(g) = \frac{1}{2} \|u_{h\alpha g} - z_d\|_H^2 + \frac{M}{2} \|g\|_H^2 \qquad (18)$$

where $u_{hg}$ and $u_{h\alpha g}$ are the discrete system states defined as the solution of the following discrete elliptic variational equalities [16, 24]:

$$u_{hg} \in K_h : \quad a(u_{hg}, v_h) = (g, v_h)_H - \int_{\Gamma_2} q v_h d\gamma, \quad \forall v_h \in V_{0h}, \qquad (19)$$

$$u_{h\alpha g} \in V_h : \quad a_\alpha(u_{h\alpha g}, v_h) = (g, v_h)_H - \int_{\Gamma_2} q v_h d\gamma + \alpha \int_{\Gamma_1} b v_h d\gamma, \quad \forall v_h \in V_h. \qquad (20)$$

The corresponding discrete distributed optimal control problems consists in finding $g_{h_{op}}, g_{h\alpha_{op}} \in H$ such that:

$$\text{Problem } (P_h): \quad J_h\left(g_{h_{op}}\right) = \underset{g \in H}{\text{Min}} J_h(g), \quad \text{Problem } (P_{h\alpha}): \quad J_{h\alpha}\left(g_{h\alpha_{op}}\right) = \underset{g \in H}{\text{Min}} J_{h\alpha}(g). \qquad (21)$$

and their corresponding discrete adjoint states $p_{hg}$ and $p_{h\alpha g}$ are defined respectively as the solution of the following discrete elliptic variational equalities:

$$p_{hg} \in V_{0h}: \quad a(p_{hg}, v_h) = (u_{hg} - z_d, v_h), \quad \forall v_h \in V_{0h} \qquad (22)$$

$$p_{h\alpha g} \in V_h: \quad a_\alpha(p_{h\alpha g}, v_h) = (u_{h\alpha g} - z_d, v_h), \quad \forall v_h \in V_h. \qquad (23)$$

**Remark 1**. We note that the discrete (in the n-dimensional space) distributed optimal control problem ($P_h$) and ($P_{h\alpha}$) are still an infinite dimensional optimal control problem since the control space is not discretized.

**Lemma 2**
*(i) There exist unique solutions $u_{hg} \in K_h$ and $p_{hg} \in V_{0h}$ of the elliptic variational equalities (19) and (22) respectively, and $u_{h\alpha g} \in V_h$ and $p_{h\alpha g} \in V_h$ of the elliptic variational equalities (20), and (23) respectively $\forall g \in H$, $\forall q \in Q$, $b > 0$ on $\Gamma_1$.*
*(ii) The operators $g \in H \to u_{hg} \in V$, and $g \in H \to u_{h\alpha g} \in V$ are Lipschitzians, i.e.*

$$a) \ \|u_{hg_2} - u_{hg_1}\|_v \leq \frac{1}{\lambda} \|g_2 - g_1\|_H, \quad b) \ \|u_{h\alpha g_2} - u_{h\alpha g_1}\|_v \leq \frac{1}{\lambda_\alpha} \|g_2 - g_1\|_H, \quad \forall g_1, g_2 \in H, \forall h > 0. \qquad (24)$$

*(iii) The operators $g \in H \to p_{hg} \in V_{0g}$, and $g \in H \to p_{h\alpha g} \in V_h$ are Lipschitzians and strictly monotones, i.e.*

$$\left(p_{hg_2} - p_{hg_1}, g_2 - g_1\right)_H = \|u_{hg_2} - u_{hg_1}\|_H^2 \geq 0, \quad \left(p_{h\alpha g_2} - p_{h\alpha g_1}, g_2 - g_1\right)_H = \|u_{h\alpha g_2} - u_{h\alpha g_1}\|_H^2 \geq 0, \quad \forall g_1, g_2 \in H, \forall h > 0 \qquad (25)$$

$$a) \ \|p_{hg_2} - p_{hg_1}\|_V \leq \frac{1}{\lambda^2}\|g_2 - g_1\|_H, \quad b) \ \|p_{h\alpha g_2} - p_{h\alpha g_1}\|_V \leq \frac{1}{\lambda_\alpha^2}\|g_2 - g_1\|_H, \ \forall g_1, g_2 \in H, \forall h > 0. \quad (26)$$

**Proof.** We use the Lax-Milgram Theorem, the variational equalities (19), (20), (22) and (23), the coerciveness (11) and (12) and following [10, 18, 25]. □

**Theorem 3**

*(i) The discrete cost functional $J_h$ and $J_{h\alpha}$ are H - elliptic and strictly convexe applications, that is:*

$$(1-t)J_h(g_2) + tJ_h(g_1) - J_h(tg_1 + (1-t)g_2) = \frac{t(1-t)}{2}\|u_{hg_2} - u_{hg_1}\|_H^2 + M\frac{t(1-t)}{2}\|g_2 - g_1\|_H^2$$
$$\geq M\frac{t(1-t)}{2}\|g_2 - g_1\|_H^2, \quad \forall g_1, g_2 \in H, \forall t \in [0,1]. \quad (27)$$

$$(1-t)J_{h\alpha}(g_2) + tJ_{h\alpha}(g_1) - J_{h\alpha}(tg_1 + (1-t)g_2) = \frac{t(1-t)}{2}\|u_{h\alpha g_2} - u_{h\alpha g_1}\|_H^2 + M\frac{t(1-t)}{2}\|g_2 - g_1\|_H^2$$
$$\geq M\frac{t(1-t)}{2}\|g_2 - g_1\|_H^2, \quad \forall g_1, g_2 \in H, \forall t \in [0,1]. \quad (28)$$

*(ii) There exist a unique optimal controls $g_{h_{op}} \in H$ and $g_{h\alpha_{op}} \in H$ that satisfy the optimization problems (21a) and (21b) respectively.*

*(iii) $J_h$ and $J_{h\alpha}$ are Gâteaux differenciable applications and their derivatives are given by the following expressions:*

$$a) \ J_h'(g) = Mg + p_{hg}, \quad b) \ J_{h\alpha}'(g) = Mg + p_{h\alpha g}, \quad \forall g \in H, \quad \forall h > 0 \quad (29)$$

*(iv) The optimality condition for the optimization problems (21a) and (21b) are given by:*

$$a) \ J_h'(g_{h_{op}}) = 0 \Leftrightarrow g_{h_{op}} = -\frac{1}{M}p_{hg_{h_{op}}}, \quad b) \ J_{h\alpha}'(g_{h\alpha_{op}}) = 0 \Leftrightarrow g_{h\alpha_{op}} = -\frac{1}{M}p_{h\alpha g_{h\alpha_{op}}}. \quad (30)$$

*(v) $J_h'$ and $J_{h\alpha}'$ are Lipschitzians and strictly monotone operators, i.e.*

$$\|J_h'(g_2) - J_h'(g_1)\|_H \leq \left(M + \frac{1}{\lambda^2}\right)\|g_2 - g_1\|_H, \quad \forall g_1, g_2 \in H, \forall h > 0 \quad (31)$$

$$\langle J_h'(g_2) - J_h'(g_1), g_2 - g_1 \rangle = \|u_{hg_2} - u_{hg_1}\|_H^2 + M\|g_2 - g_1\|_H^2 \geq M\|g_2 - g_1\|_H^2, \forall g_1, g_2 \in H, \forall h > 0 \quad (32)$$

$$\|J_{h\alpha}'(g_2) - J_{h\alpha}'(g_1)\|_H \leq \left(M + \frac{1}{\lambda_\alpha^2}\right)\|g_2 - g_1\|_H, \quad \forall g_1, g_2 \in H, \forall h > 0 \quad (33)$$

$$\langle J_{h\alpha}'(g_2) - J_{h\alpha}'(g_1), g_2 - g_1 \rangle = \|u_{h\alpha g_2} - u_{h\alpha g_1}\|_H^2 + M\|g_2 - g_1\|_H^2 \geq M\|g_2 - g_1\|_H^2, \forall g_1, g_2 \in H, \forall h > 0. \quad (34)$$

**Proof.** We use the definitions (18a,b), the elliptic variational equalities (19) and (20) and the coerciveness (11) and (12), following [10,18,25]. Moreover, the functional $J_h'$ and $J_{h\alpha}'$ are given by:

$$\langle J_h'(g), f \rangle = \lim_{t \to 0^+} \frac{J_h(g + tf) - J_h(g)}{t} = (Mg + p_{hg}, f)_H, \quad \forall g, f \in H, \quad (35)$$

$$\langle J_{h\alpha}'(g), f \rangle = \lim_{t \to 0^+} \frac{J_{h\alpha}(g + tf) - J_{h\alpha}(g)}{t} = (Mg + p_{h\alpha g}, f)_H, \quad \forall g, f \in H. \ \square \quad (36)$$

We define the operators:

$$W_h : H \to V_{0h} \subset V_0 \subset H / W_h(g) = -\frac{1}{M}p_{hg}, \quad (37)$$

$$W_{h\alpha} : H \to V_h \subset V \subset H / W_{h\alpha}(g) = -\frac{1}{M}p_{h\alpha g}. \quad (38)$$

**Theorem 4**

*We have that:*

*(i) $W_h$ and $W_{h\alpha}$ are Lipschitzian operators, that is:*

$$a) \ \|W_h(g_2) - W_h(g_1)\|_V \leq \frac{1}{M\lambda^2}\|g_2 - g_1\|_H, \quad b) \ \|W_{h\alpha}(g_2) - W_{h\alpha}(g_1)\|_V \leq \frac{1}{M\lambda_\alpha^2}\|g_2 - g_1\|_H, \forall g_1, g_2 \in H, \forall h > 0 \quad (39)$$

*(ii) $W_h$ ($W_{h\alpha}$) is a contraction operator if and only if $M$ is large, that is:*

$$\text{a) } M > \frac{1}{\lambda^2}, \quad \text{b) } M > \frac{1}{\lambda_\alpha^2}. \tag{40}$$

*(iii) If M verifies the inequality (40) then the discrete distributional optimal control $g_{h_{op}} \in H$ ( $g_{h\alpha_{op}} \in H$ ) can be also obtained as the unique fixed point of the operator $W_h$ ($W_{h\alpha}$), that is:*

$$g_{h_{op}} = -\frac{1}{M} p_{hg_{h_{op}}} \Leftrightarrow W_h\left(g_{h_{op}}\right) = g_{h_{op}}, \quad g_{h\alpha_{op}} = -\frac{1}{M} p_{h\alpha g_{h\alpha_{op}}} \Leftrightarrow W_{h\alpha}\left(g_{h\alpha_{op}}\right) = g_{h\alpha_{op}} \tag{41}$$

**Proof.** We use the definitions (37) and (38), and the properties (26a,b) and (41a,b). □

## III. Convergence of the Discrete Distributed Optimal Control Problems $(P_h)$ to $(P)$ and $(P_{h\alpha})$ to $(P_\alpha)$ when $h \to 0$

We obtain the following error estimations between the continuous and discrete solutions:

**Lemma 5**
*(i) If the continuous system states and the continuous adjoint system states have the regularity $u_g, u_{\alpha g} \in H^r(\Omega), p_g, p_{\alpha g} \in H^r(\Omega) (1 < r \leq 2)$ then $\forall \alpha > 0, \forall g \in H$ we have the following estimations:*

$$\|u_g - u_{hg}\|_V \leq \frac{c_0}{\sqrt{\lambda}} \|u_g\|_r h^{r-1}, \quad \|p_g - p_{hg}\|_V \leq ch^{r-1}, \quad \forall g \in H, h > 0, \tag{42}$$

$$\|u_{h\alpha g} - u_{\alpha g}\|_V \leq ch^{r-1}, \quad \|p_{h\alpha g} - p_{\alpha g}\|_V \leq ch^{r-1}, \tag{43}$$

*where c's are constants independents of h.*
*(ii) We have the following convergences:*

$$\lim_{h \to 0^+} \|u_g - u_{hg}\|_V = 0, \quad \lim_{h \to 0^+} \|p_g - p_{hg}\|_V = 0, \quad \forall g \in H, \tag{44}$$

$$\lim_{h \to 0^+} \|u_{h\alpha g} - u_{\alpha g}\|_V = 0, \quad \lim_{h \to 0^+} \|p_{h\alpha g} - p_{\alpha g}\|_V = 0, \quad \forall \alpha > 0, \forall g \in H. \tag{45}$$

**Proof.** We use the variational equalities (5), (7) and (19), $v_h = \pi_h(u_g)$ in the variational equalities (19) and (20), the coerciveness properties (11) and (12), the estimations (21) and we have the following properties:

$$a\left(p_g - p_{hg}, \pi_h(p_g) - p_{hg}\right) = \left(u_g - u_{hg}, \pi_h(p_g) - p_{hg}\right), \quad \forall g \in H, \tag{46}$$

$$a_\alpha\left(p_{\alpha g} - p_{h\alpha g}, \pi_h(p_{\alpha g}) - p_{h\alpha g}\right) = \left(u_{h\alpha g} - u_{\alpha g}, \pi_h(p_{\alpha g}) - p_{h\alpha g}\right)_H, \quad \forall g \in H, \tag{47}$$

following a similar method given in [25] for Neumann boundary optimal control problems. □

**Theorem 6**
We consider the continuous system states and adjoint system states have the regularities $u_g, u_{\alpha g_{\alpha_{op}}} \in H^r(\Omega)$ and $p_g, p_{\alpha g_{\alpha_{op}}} \in H^r(\Omega)$ $(1 < r \leq 2)$.

*i) We have the following limits:*

$$\lim_{h \to 0^+} \|g_{h_{op}} - g_{op}\|_V = 0, \quad \lim_{h \to 0^+} \|u_{hg_{h_{op}}} - u_{g_{op}}\|_V = 0, \quad \lim_{h \to 0^+} \|p_{hg_{h_{op}}} - p_{g_{op}}\|_V = 0, \tag{48}$$

$$\lim_{h \to 0^+} \|g_{h\alpha_{op}} - g_{\alpha_{op}}\|_H = 0, \quad \lim_{h \to 0^+} \|u_{h\alpha g_{h\alpha_{op}}} - u_{\alpha g_{\alpha_{op}}}\|_V = 0, \quad \lim_{h \to 0^+} \|p_{h\alpha g_{h\alpha_{op}}} - p_{\alpha g_{\alpha_{op}}}\|_V = 0, \forall \alpha > 1 \tag{49}$$

*ii) If M verifies the inequalities (40a,b) then we have the following error bonds:*

$$\|g_{h_{op}} - g_{op}\|_H \leq ch^{r-1}, \|u_{hg_{h_{op}}} - u_{g_{op}}\|_V \leq ch^{r-1}, \|p_{hg_{h_{op}}} - p_{g_{op}}\|_V \leq ch^{r-1} \tag{50}$$

$$\|g_{h\alpha_{op}} - g_{\alpha_{op}}\|_H \leq ch^{r-1}, \|u_{h\alpha g_{h\alpha_{op}}} - u_{\alpha g_{\alpha_{op}}}\|_V \leq ch^{r-1}, \|p_{h\alpha g_{h\alpha_{op}}} - p_{\alpha g_{\alpha_{op}}}\|_V \leq ch^{r-1} \tag{51}$$

*where c's are constants independents of h.*

**Proof.** We follow a similar method to the one developed in [25] for Neumann boundary optimal control problems by using the elliptic variational equalities (19), (20), (22) and (23), the thesis holds. □

**Remark 2.** If M verifies the inequalities (40a,b) we can obtain (48) and (49) by using the characterization of the fixed point (41a,b), and the uniqueness of the optimal controls $g_{op} \in H$ and $g_{\alpha_{op}} \in H$.

**Remark 3**. Firstly, In Theorem 6 we have used the restriction $\alpha > 1$ by splitting the bilinear form $a_\alpha$ by [21, 24,25]

$$a_\alpha(u,v) = a_1(u,v) + (\alpha - 1)\int_{\Gamma_1} uv\, d\gamma, \tag{52}$$

then it can be replaced by $\alpha \geq \alpha_0$ for any $\alpha_0 > 0$.

Now, we have some estimations for the discrete cost functional $J_{h\alpha}$ and $J_h$ when $h \to 0$.

**Lemma 7** *If M verifies the inequality (40a,b) and the continuous system states and adjoint system states have the regularities $u_g, u_{\alpha g} \in H^r(\Omega), p_g, p_{\alpha g} \in H^r(\Omega) (1 < r \leq 2)$ then we have the following error bonds:*

$$\frac{M}{2}\left\|g_{h_{op}} - g_{op}\right\|_H^2 \leq J(g_{h_{op}}) - J(g_{op}) \leq Ch^{2(r-1)}, \quad \frac{M}{2}\left\|g_{h\alpha_{op}} - g_{\alpha_{op}}\right\|_H^2 \leq J_\alpha(g_{h\alpha_{op}}) - J_\alpha(g_{\alpha_{op}}) \leq Ch^{2(r-1)} \tag{53}$$

$$\frac{M}{2}\left\|g_{h_{op}} - g_{op}\right\|_H^2 \leq J_h(g_{op}) - J_h(g_{h_{op}}) \leq Ch^{2(r-1)}, \quad \frac{M}{2}\left\|g_{h\alpha_{op}} - g_{\alpha_{op}}\right\|_H^2 \leq J_{h\alpha}(g_{\alpha_{op}}) - J_{h\alpha}(g_{h\alpha_{op}}) \leq Ch^{2(r-1)} \tag{54}$$

$$\left|J_h(g_{op}) - J(g_{op})\right| \leq Ch^{r-1}, \quad \left|J_h(g_{h_{op}}) - J(g_{op})\right| \leq Ch^{r-1} \tag{55}$$

$$\left|J_{h\alpha}(g_{op}) - J_\alpha(g_{op})\right| \leq Ch^{r-1}, \quad \left|J_{h\alpha}(g_{h\alpha_{op}}) - J_\alpha(g_{\alpha_{op}})\right| \leq Ch^{r-1} \tag{56}$$

*where C's are constants independents of h and $\alpha$.*

**Proof.** Estimations (53) and (54) follow from the estimations (42) and (50), and the equalities:

$$J(g_{h_{op}}) - J(g_{op}) = \frac{1}{2}\left\|u_{g_{h_{op}}} - u_{g_{op}}\right\|_H^2 + \frac{M}{2}\left\|g_{h_{op}} - g_{op}\right\|_H^2, \tag{57}$$

$$J_h(g_{op}) - J_h(g_{h_{op}}) = \frac{1}{2}\left\|u_{hg_{h_{op}}} - u_{hg_{op}}\right\|_H^2 + \frac{M}{2}\left\|g_{h_{op}} - g_{op}\right\|_H^2. \tag{58}$$

$$J_\alpha(g_{h\alpha_{op}}) - J(g_{\alpha_{op}}) = \frac{1}{2}\left\|u_{\alpha g_{h\alpha_{op}}} - u_{\alpha g_{op}}\right\|_H^2 + \frac{M}{2}\left\|g_{h\alpha_{op}} - g_{\alpha_{op}}\right\|_H^2, \tag{59}$$

$$J_{h\alpha}(g_{\alpha_{op}}) - J_{h\alpha}(g_{h\alpha_{op}}) = \frac{1}{2}\left\|u_{h\alpha g_{h_{op}}} - u_{h\alpha g_{h\alpha_{op}}}\right\|_H^2 + \frac{M}{2}\left\|g_{h\alpha_{op}} - g_{\alpha_{op}}\right\|_H^2. \tag{60}$$

Estimations (55) follow from the estimations (24), (42) and (50), the triangular inequality for norms and the inequalities:

$$\left|J_h(g) - J(g)\right| \leq \left(\frac{1}{2}\left\|u_{hg} - u_g\right\|_H + \left\|u_g - z_d\right\|_H\right)\left\|u_{hg} - u_g\right\|_H, \quad \forall g \in H, \tag{61}$$

$$\left|J_{h\alpha}(g) - J_\alpha(g)\right| \leq \left(\frac{1}{2}\left\|u_{h\alpha g} - u_{\alpha g}\right\|_H + \left\|u_{\alpha g} - z_d\right\|_H\right)\left\|u_{h\alpha g} - u_{\alpha g}\right\|_H, \quad \forall g \in H. \,\square \tag{62}$$

## IV. Convergence of the Discrete Optimal Control Problems $(P_{h\alpha})$ to $(P_h)$ when $\alpha \to +\infty$

For a fixed $h > 0$ we have:

**Lemma 8.** For a fixed $g \in H$ we have the following limits:

$$a)\ \lim_{\alpha \to +\infty}\left\|u_{h\alpha g} - u_{hg}\right\|_V = 0, \quad b)\ \lim_{\alpha \to +\infty}\left\|p_{h\alpha g} - p_{hg}\right\|_V = 0, \quad \forall g \in H, \forall h > 0, \tag{63}$$

**Proof.** For fixed $g \in H, h > 0$, and by using the variational equalities (19) and (20), and the splitting form (52) we obtain the following estimations:

$$\left\|u_{h\alpha g} - u_{hg}\right\|_V \leq c, \quad (\alpha - 1)\int_{\Gamma_1}(u_{h\alpha g} - b)^2 d\gamma \leq c, \quad \forall \alpha > 1. \tag{64}$$

From the above inequalities we deduce that:

$$\exists \eta_{hg} \in V / u_{h\alpha g} \longrightarrow \eta_{hg} \text{ in } V \text{ weak (in } H \text{ strong) as } \alpha \to +\infty \text{ with } \eta_{hg}/\Gamma_1 = b. \tag{65}$$

By using the variational equality (20) we can pass to the limit when $\alpha \to +\infty$, and by uniqueness of the variational equality (19) we obtain that $\eta_{hg} = u_{hg}$. By using the above properties, and the variational equalities (19) and (20), we deduce (63a) and by using a similar method we can obtain the limit (63b) for the discrete adjoint system state. □

**Theorem 9**

We have the following limits:

$$\lim_{\alpha \to +\infty} \left\| u_{h\alpha g_{h\alpha_{op}}} - u_{hg_{h_{op}}} \right\|_V = \lim_{\alpha \to +\infty} \left\| p_{h\alpha g_{h\alpha_{op}}} - p_{hg_{h_{op}}} \right\|_V = \lim_{\alpha \to +\infty} \left\| g_{h\alpha_{op}} - g_{h_{op}} \right\|_H = 0, \quad \forall h > 0. \tag{66}$$

**Proof.** We omit this proof because we prefer to prove the next one with more details.

## V Double Convergence of the Discrete Distributed Optimal Control Problem $(P_{h\alpha})$ to $(P)$ when $(h, \alpha) \to (0, +\infty)$

For the discrete distributed optimal control problem $(P_{h\alpha})$ we will now consider the double limit $(h, \alpha) \to (0, +\infty)$.

**Theorem 10**

We have the following limits:

$$\lim_{(h,\alpha) \to (0,+\infty)} \left\| u_{h\alpha g_{h\alpha_{op}}} - u_{g_{op}} \right\|_V = \lim_{(h,\alpha) \to (0,+\infty)} \left\| p_{h\alpha g_{h\alpha_{op}}} - p_{g_{op}} \right\|_V = \lim_{(h,\alpha) \to (0,+\infty)} \left\| g_{h\alpha_{op}} - g_{op} \right\|_H = 0. \tag{67}$$

**Proof.** From now on we consider that c's represent positive constants independents simultaneously of $h > 0$ and $\alpha > 0$. We note that without loss of generality we can consider $\alpha > 1$ by splitting suitably the bilinear form $a_\alpha$ taking into account (52). We show a sketch of the proof by obtaining the following estimations (for $\forall h > 0$ and $\forall \alpha > 1$):

$$\|u_{h0}\|_V \leq c_1 = b + \frac{\|q\|_Q \|\gamma_0\|}{\lambda}, \tag{68}$$

$$\|u_{h\alpha 0}\|_V \leq c_2 = b\left(1 + \frac{1}{\lambda_1}\right) + \|q\|_Q \|\gamma_0\| \left[\frac{1}{\lambda_1} + \frac{1}{\lambda}\left(1 + \frac{1}{\lambda_1}\right)\right], \tag{69}$$

$$(\alpha - 1) \int_{\Gamma_1} (u_{h\alpha 0} - b)^2 \, d\gamma \leq c_3 = \frac{\lambda_1}{\lambda^2} \left[ b + \|q\|_Q \|\gamma_0\| \left(1 + \frac{1}{\lambda}\right) \right]^2, \tag{70}$$

$$\|g_{h\alpha_{op}}\|_H \leq c_4 = \sqrt{\frac{1}{M}} \left( \|z_d\| + c_2 \right), \tag{71}$$

$$\|u_{h\alpha g_{h\alpha_{op}}}\|_H \leq c_5 = 2\|z_d\| + c_2, \tag{72}$$

$$\|g_{h_{op}}\|_H \leq c_6 = \sqrt{\frac{1}{M}} \left( \|z_d\| + c_1 \right), \tag{73}$$

$$\|u_{hg_{h_{op}}}\|_V \leq c_7 = b\left(1 + \frac{1}{\lambda\sqrt{M}}\right) + \frac{1}{\lambda\sqrt{M}} \|z_d\|_H + \frac{\|q\|_Q \|\gamma_0\|}{\lambda} \left[1 + \frac{1}{\lambda\sqrt{M}}\right], \tag{74}$$

$$\|u_{h\alpha g_{h\alpha_{op}}}\|_V \leq c_8 = \frac{\|z_d\|_H}{\sqrt{M}} \left[\frac{1}{\lambda_1} + \frac{1}{\lambda}\left(1 + \frac{1}{\lambda_1}\right)\right] + b\left(1 + \frac{1}{\lambda_1}\right)\left(1 + \frac{1}{\lambda\sqrt{M}} + \frac{1}{\lambda_1\sqrt{M}}\right)$$
$$+ \|q\|_Q \|\gamma_0\| \left[\frac{1}{\lambda_1} + \frac{1}{\lambda}\left(1 + \frac{1}{\lambda_1}\right) + \frac{1}{\sqrt{M}}\left(\frac{1}{\lambda_1^2} + \frac{1}{\lambda^2}\left(1 + \frac{1}{\lambda_1}\right) + \frac{1}{\lambda\lambda_1}\left(1 + \frac{1}{\lambda_1}\right)\right)\right], \tag{75}$$

$$(\alpha - 1) \int_{\Gamma_1} (u_{h\alpha g_{h\alpha_{op}}} - b)^2 \, d\gamma \leq c_9 = \frac{1}{\lambda_1} \left[ \frac{\|z_d\|_H}{\sqrt{M}}\left(1 + \frac{1}{\lambda}\right) + b\left(1 + \frac{1}{\sqrt{M}}\left(1 + \frac{1}{\lambda_1} + \frac{1}{\lambda}\right)\right) \right.$$
$$\left. + \|q\|_Q \|\gamma_0\| \left(1 + \frac{1}{\lambda} + \frac{1}{\sqrt{M}}\left(\frac{1}{\lambda} + \frac{1}{\lambda_1} + \frac{1}{\lambda^2} + \frac{1}{\lambda\lambda_1}\right)\right) \right]^2, \tag{76}$$

$$\left\| p_{hg_{h_{op}}} \right\|_V \leq c_{10} = \frac{1}{\lambda}\left(1+\frac{1}{\lambda\sqrt{M}}\right)\left(\|z_d\|_H + c_1\right). \tag{77}$$

$$\left\| p_{h\alpha g_{h\alpha_{op}}} \right\|_V \leq c_{11} = \|z_d\|_H \left[\frac{1}{\lambda_1}\left(1+\frac{1}{\sqrt{M}}\left(\frac{1}{\lambda_1}+\frac{1}{\lambda}+\frac{1}{\lambda\lambda_1}\right)\right)+\frac{1}{\lambda}\left(1+\frac{1}{\lambda_1}\right)\left(1+\frac{1}{\lambda\sqrt{M}}\right)\right]$$
$$+b\left[\frac{1}{\lambda_1}\left(1+\frac{1}{\lambda_1}\right)\left(1+\frac{1}{\lambda\sqrt{M}}+\frac{1}{\lambda_1\sqrt{M}}\right)+\frac{1}{\lambda}\left(1+\frac{1}{\lambda_1}\right)\left(1+\frac{1}{\lambda\sqrt{M}}\right)\right]$$
$$+\|q\|_Q \|\gamma_0\|\left[\frac{1}{\lambda_1}\left[\frac{1}{\lambda_1}+\frac{1}{\lambda}\left(1+\frac{1}{\lambda_1}\right)+\frac{1}{\sqrt{M}}\left(\frac{1}{\lambda\lambda_1}+\frac{1}{\lambda^2}\left(1+\frac{1}{\lambda_1}\right)+\frac{1}{\lambda_1^2}\left(1+\frac{1}{\lambda}\right)\right)\right]+\frac{1}{\lambda^2}\left(1+\frac{1}{\lambda_1}\right)\left(1+\frac{1}{\lambda\sqrt{M}}\right)\right]$$
$$\tag{78}$$

Therefore, from the above estimations we have that:

$$\exists f \in H \ / \ g_{h\alpha_{op}} \longrightarrow f \text{ in } H \text{ weak as } (h,\alpha) \to (0,+\infty), \tag{79}$$

$$\exists \eta \in V \ / \ u_{h\alpha g_{h\alpha_{op}}} \longrightarrow \eta \text{ in } V \text{ weak (in } H \text{ strong) as } (h,\alpha) \to (0,+\infty) \text{ with } \eta/\Gamma_1 = b, \tag{80}$$

$$\exists \xi \in V \ / \ p_{h\alpha g_{h\alpha_{op}}} \longrightarrow \xi \text{ in } V \text{ weak (in } H \text{ strong) as } (h,\alpha) \to (0,+\infty) \text{ with } \xi/\Gamma_1 = 0, \tag{81}$$

and

$$\exists f_h \in H \ / \ g_{h\alpha_{op}} \longrightarrow f_h \text{ in } H \text{ weak as } \alpha \to +\infty, \tag{82}$$

$$\exists \eta_h \in V \ / \ u_{h\alpha g_{h\alpha_{op}}} \longrightarrow \eta_h \text{ in } V \text{ weak (in } H \text{ strong) as } \alpha \to +\infty \text{ with } \eta_h/\Gamma_1 = b, \tag{83}$$

$$\exists \xi_h \in V \ / \ p_{h\alpha g_{h\alpha_{op}}} \longrightarrow \xi_h \text{ in } V \text{ weak (in } H \text{ strong) as } \alpha \to +\infty \text{ with } \xi_h/\Gamma_1 = 0, \tag{84}$$

and

$$\exists f_\alpha \in H \ / \ g_{h\alpha_{op}} \longrightarrow f_\alpha \text{ in } H \text{ weak as } h \to 0, \tag{85}$$

$$\exists \eta_\alpha \in V \ / \ u_{h\alpha g_{h\alpha_{op}}} \longrightarrow \eta_\alpha \text{ in } V \text{ weak (in } H \text{ strong) as } h \to 0 \text{ with } \eta_\alpha/\Gamma_1 = b, \tag{86}$$

$$\exists \xi_\alpha \in V \ / \ p_{h\alpha g_{h\alpha_{op}}} \longrightarrow \xi_\alpha \text{ in } V \text{ weak (in } H \text{ strong) as } h \to 0 \text{ with } \xi_\alpha/\Gamma_1 = 0. \tag{87}$$

Taking into account the uniqueness of the distributed optimal control problems $(P_{h\alpha}), (P_\alpha), (P_h)$ and $(P)$, and the uniqueness of the elliptic variational equalities corresponding to their state systems we get

$$\eta_h = u_{hf_h} = u_{hg_{h_{op}}}, \quad \xi_h = p_{hf_h} = p_{hg_{h_{op}}}, \quad f_h = g_{h_{op}}, \tag{88}$$

$$\eta_\alpha = u_{\alpha f_\alpha} = u_{\alpha g_{\alpha_{op}}}, \quad \xi_\alpha = p_{\alpha f_\alpha} = p_{\alpha g_{\alpha_{op}}}, \quad f_\alpha = g_{\alpha_{op}}, \tag{89}$$

and the limits (48). Now, by using [11] we get

$$\lim_{\alpha \to +\infty}\left\|f_\alpha - g_{op}\right\|_H = 0, \quad \lim_{\alpha \to +\infty}\left\|\eta_\alpha - u_{g_{op}}\right\|_V = 0, \quad \lim_{\alpha \to +\infty}\left\|\xi_\alpha - p_{g_{op}}\right\|_V = 0, \tag{90}$$

and therefore the three double limits (67) hold when $(h,\alpha) \to (0,+\infty)$.

**Remark 4.** We note that this double convergence is a novelty with respect to the recent results obtained for a family of discrete Neumann boundary optimal control problems [25].

**Proof of Theorem 1.** It is a consequence of the properties (48), (49), (66), (67) and [10,11].

## Conclusions

We have studied the numerical analysis of the discrete distributed optimal control problems $(P_h)$ and $(P_{h\alpha})$, and the corresponding asymptotic behaviour when $\alpha \to \infty$, $h \to 0$ and $(h,\alpha) \to (0,+\infty)$ by using the finite element method. We have defined the discrete cost functional $J_h$ and $J_{h\alpha}$, the discrete variational equalities for the system states $u_{hg}$ and $u_{h\alpha g}$ for each $\alpha, h > 0$ and $g \in H$, and the discrete variational equalities for the adjoint system states $p_{hg}$ and $p_{h\alpha g}$ for each $\alpha, h > 0$, and $g \in H$. We have characterized the discrete distributed optimal control energy $g_{h_{op}}$ and $g_{h\alpha_{op}}$ as a fixed point on $H$ of suitable discrete operators $W_h$ and $W_{h\alpha}$ over his adjoint system states $p_{hg_{op}}$ and $p_{h\alpha g_{h\alpha_{op}}}$ respectively for each $\alpha > 0$. We have also studied the convergence of the discrete distributed optimal control problems ($P_{h\alpha}$) to ($P_h$) when $\alpha \to \infty$ for each $h > 0$, the convergence of the discrete distributed optimal control problems ($P_{h\alpha}$) to ($P_\alpha$), and ($P_h$) to $P$ when $h \to 0$ for each $\alpha > 0$, and the double convergence of the discrete distributed optimal control problems ($P_{h\alpha}$) to ($P$) when $(h, \alpha) \to (0, +\infty)$. Thus, we have obtained a commutative diagram (see Introduction) which relates the continuous and discrete distributed mixed optimal control problems $(P_{h\alpha}), (P_\alpha), (P_h)$ and ($P$) by taking the limits $h \to 0$, $\alpha \to \infty$ and $(h, \alpha) \to (0, +\infty)$ respectively.

## Acknowledgements

The present work has been partially sponsored by the Projects PIP No 0534 from CONICET - Univ. Austral, Rosario, Argentina, and AFOSR-SOARD Grant FA9550-14-1-0122.

## References

[1] A. Azzam, E. Kreyszig, "On solutions of elliptic equations satisfying mixed boundary conditions", SIAM J. Math. Anal., 13 (1982), 254-262.
[2] M. Bergounioux, "Optimal control of an obstacle problem", Appl. Math. Optim., 36 (1997), 147-172.
[3] S. Brenner, L.R. Scott, "The mathematical theory of finite element methods", Springer, New York, 2008.
[4] E. Casas, M. Mateos, "Uniform convergence of the FEM. Applications to state constrained control problems", Comput. Appl. Math., 21 (2002), 67-100.
[5] E. Casas, M. Mateos, "Dirichlet contol problems in smooth and nonsmooth convex plain domains", Control Cybernetics, 40 (2011), 931-955.
[6] E. Casas, J.P. Raymond, "Error estimates for the numerical approximation of Dirichlet boundary control for semilinear elliptic equations", SIAM J. Control Optim., 45 (2006), 1586-1611.
[7] P.G. Ciarlet, "The finite element method for elliptic problems", SIAM, Philadelphia, 2002.
[8] K. Deckelnick, A. Günther, M. Hinze, "Finite element approximation of ellliptic control problems with constraints on the gradient", Numer. Math., 111 (2009), 335-350.
[9] K. Deckelnick, M. Hinze, "Convergence of a finite element approximation to a state-constrained ellliptic control problem", SIAM J. Numer. Anal., 45 (2007), 1937-1953.
[10] C.M. Gariboldi, D.A. Tarzia, "Convergence of distributed optimal controls on the internal energy in mixed elliptic problems when the heat transfer coefficient goes to infinity", Appl. Math. Optim., 47 (2003), 213-230.
[11] C.M. Gariboldi, D.A. Tarzia, "A new proof of the convergence of the distributed optimal controls on the internal energy in mixed elliptic problems", MAT – Serie A, 7 (2004), 31-42.
[12] R. Haller-Dintelmann, C. Meyer, J. Rehberg, A. Schiela, "Hölder continuity and optimal control for nonsmooth elliptic problems", Appl. Math. Optim., 60 (2009), 397-428.
[13] M. Hintermüller, M. Hinze, "Moreau-Yosida regularization in state constrained ellliptic control problems: Error estimates and parameter adjustement", SIAM J. Numer. Anal., 47 (2009), 1666-1683.
[14] M. Hinze, "A variational discretization concept in control constrained optimization: The linear-quadratic case", Comput. Optim. Appl., 30 (2005), 45-61.
[15] M. Hinze, U. Matthes, "A note on variational dicretization of elliptic Nuemann boundary control", Control Cybernetics, 38 (2009), 577-591.
[16] D. Kinderlehrer, G. Stampacchia, "An introduction to variational inequalities and their applications", SIAM, Philadelphia, 2000.
[17] L. Lanzani, L. Capogna, R.M. Brown, "The mixed problem in $L^p$ for some two-dimensional Lipschitz domain", Math. Annalen, 342 (2008), 91-124.
[18] J.L. Lions, "Contrôle optimal des systèmes gouvernés par des équations aux dérivées partielles", Dunod, Paris, 1968.
[19] E.B. Mermri, W. Han, "Numerical approximation of a unilateral obstacle problem", J. Optim. Th. Appl., 153 (2012), 177-194.
[20] E. Shamir, "Regularization of mixed second order elliptic problems", Israel J. Math., 6 (1968), 150-168.
[21] E.D. Tabacman, D.A. Tarzia, "Sufficient and/or necessary condition for the heat transfer coefficient on $\Gamma_1$ and the heat flux on $\Gamma_2$ to obtain a steady-state two-phase Stefan problem", J. Diff. Eq., 77 (1989), 16-37.
[22] D.A. Tarzia, "An inequality for the constant heat flux to obtain a steady-state two-phase Stefan problem", Eng. Anal., 5 (1988), 177-181.
[23] D.A. Tarzia, "Numerical analysis for the heat flux in a mixed elliptic problem to obtain a discrete steady-state two-

phase Stefan problem", SIAM J. Numer. Anal., 33 (1996), 1257-1265.

[24] D.A. Tarzia, "Numerical analysis of a mixed elliptic problem with flux and convective boundary conditions to obtain a discrete solution of non-constant sign", Numer. Meth. PDE, 15 (1999), 355-369.

[25] D.A. Tarzia, "A commutative diagram among discrete and continuous boundary optimal control problems", Adv. Diff. Eq. Control Processes, 14 (2014), 23-54.

[26] F. Tröltzsch, "Optimal control of partial differential equations. Theory, methods and applications", American Mathematical Society, Providence, 2010.

[27] M. Yan, L. Chang, N. Yan, "Finite element method for constrained optimal control problems governed by nonlinear elliptic PDEs", Math. Control Related Fields, 2 (2012), 183-194.

[28] Y. Ye, C.K. Chan, H.W.J. Lee, "The existence results for obstacle optimal control problems", Appl. Math. Comput., 214 (2009), 451-456.